\documentclass[paper=a4,fontsize=12pt,abstract]{scrartcl}
\pdfoutput=1 
\usepackage[ngerman,english]{babel,varioref}
\usepackage[utf8]{inputenc} 
\usepackage[T1]{fontenc} 
\usepackage{lmodern,microtype} 
\usepackage{csquotes} 

\usepackage{amssymb}
\usepackage{mathtools}
\usepackage{amsthm}
\usepackage{pgf,tikz}
\usepackage{dsfont} 

\usepackage{setspace} 

\usepackage{booktabs} 
\usepackage{paralist}

\deffootnote{1.8em}{0em}{\textsuperscript{\thefootnotemark}\ }

\makeatletter
\g@addto@macro\bfseries{\boldmath}
\makeatother 

\newtheoremstyle{mythm}
{\topsep}   
{\topsep}   
{\itshape}  
{0pt}       
{\bfseries} 
{.}         
{5pt plus 1pt minus 1pt} 
{}          

\newtheoremstyle{mydefi}
{\topsep}   
{3ex}   
{\normalfont}  
{0pt}       
{\bfseries} 
{.}         
{5pt plus 1pt minus 1pt} 
{}          

\makeatletter
\newenvironment{myproof}[1][\proofname]{\par
	\pushQED{\qed}%
	\normalfont \topsep6\p@\@plus6\p@\relax
	\trivlist
	\item[\hskip\labelsep
	\itshape
	#1\@addpunct{.}]\ignorespaces
}{%
	\popQED\endtrivlist\@endpefalse\bigskip
}
\makeatother

\theoremstyle{mythm}
\newtheorem{satz}{Satz}[section]
\newtheorem{thm}[satz]{Theorem}
\newtheorem{prop}[satz]{Proposition}

\newtheorem{rem}[satz]{Remark}

\theoremstyle{mydefi}
\newtheorem{defi}[satz]{Definition}

\newtheorem{ex}[satz]{Example}

\DeclareMathOperator{\Aut}{Aut}
\DeclareMathOperator{\SL}{SL}

\DeclareMathOperator{\PSL}{PSL}

\newcommand{\inv}[1]{#1^{-1}}
\newcommand{\mat}[4]{\left( \begin{smallmatrix}#1&#2\\#3&#4\end{smallmatrix} \right)}
\newcommand{\set}[2]{\left\{{#1}\left|\vphantom{#1#2\strut}\right.\,{#2\strut}\right\}}
\newcommand{\medset}[2]{\set{\smash{#1}}{#2}}

\usepackage[backend=biber,hyperref]{biblatex}
\usepackage[colorlinks=true,linkcolor=.,citecolor=.,urlcolor=blue,breaklinks=true,%
pdftitle={Three Affine SL(2,8)-Unitals},%
pdfauthor={Verena Möhler}%
]{hyperref}

\title{Three Affine \texorpdfstring{$\SL(2,8)$-Unitals}{SL(2,8)-Unitals}}
\author{Verena Möhler}
\date{December 17, 2020}

\begin{document}

\maketitle

\begin{abstract}
	\noindent $\SL(2,q)$-unitals are unitals of order $q$ admitting a regular action of $\SL(2,q)$ on the 
	complement of some block. We introduce three non-classical affine $\SL(2,8)$-unitals and their full automorphism groups. Each of those three affine unitals can be completed to at least two non-isomorphic unitals, leading to six pairwise non-isomorphic unitals of order $8$.\bigskip
	
	\smallskip
	
\end{abstract}

Most of the results in the present paper have been obtained in the author's Ph.\,D.\ thesis \cite{diss}, where detailed arguments can be found for some statements that
we leave to the reader here.

\section{Preliminaries}

One strategy to construct projective planes is to build an affine plane first and then to add \emph{points at infinity}, namely a new point for each parallel class and a line containing all these new points. This strategy of constructing an affine part of a geometry first and then completing it by adding some objects at infinity can successfully be applied to other incidence structures than affine and projective planes. We apply such an approach to unitals.\bigskip

A \textbf{unital of order $n$} is a $2$-$(n^3+1,n+1,1)$ design, i.\,e.\ an incidence structure with $n^3+1$ points, $n+1$ points on each block and unique joining blocks for any two points. We consider \emph{affine unitals}, which arise from unitals by removing one block (and all the points on it) and can be completed to unitals via a parallelism on the short blocks. We give an axiomatic description:

\goodbreak
\begin{defi}
	Let $n\in\mathbb N$, $n\geq2$. An incidence structure $\mathbb U=(\mathcal P,\mathcal B,\emph I)$ is called an \textbf{affine unital of order \boldmath $n$} if:
	\begin{itemize}
		\item[(AU1)] There are $n^3-n$ points.
		\item[(AU2)] Each block is incident with either $n$ or $n+1$ points. The blocks incident with $n$~points will be called \textbf{short blocks} and the blocks incident with $n+1$ points will be called \textbf{long blocks}.
		\item[(AU3)] Each point is incident with $n^2$ blocks.
		\item[(AU4)] For any two points there is exactly one block incident with both of them.
		\item[(AU5)] There exists a \textbf{parallelism} on the short blocks, meaning a partition of the set of all short blocks into $n+1$ parallel classes of size $n^2-1$ such that the blocks of each parallel class are pairwise non-intersecting.
	\end{itemize}	
\end{defi}

The existence of a parallelism as in (AU5) must explicitly be required (see \cite[Example 3.10]{diss}). An affine unital $\mathbb U$ of order $n$ with parallelism $\pi$ can be completed to a unital $\mathbb U^\pi$ of order $n$ as follows: For each parallel class, add a new point that is incident with each short block of that class. Then add a single new block $[\infty]^\pi$, incident with the $n+1$ new points (see \cite[Proposition 3.9]{diss}). We call $\mathbb U^\pi$ the \textbf{$\pi$-closure} of $\mathbb U$. Note that the closure depends on the parallelism, which need not be unique. Given an affine unital $\mathbb U$ with parallelisms $\pi$ and $\pi'$, the closures $\mathbb U^\pi$ and $\mathbb U^{\pi'}$ are isomorphic with $[\infty]^\pi \mapsto [\infty]^{\pi'}$ exactly if there is an automorphism of $\mathbb U$ which maps $\pi$ to $\pi'$ (see \cite[Proposition 3.12]{diss}).

\section{Affine \texorpdfstring{$\SL(2,q)$-Unitals}{SL(2,q)-Unitals}}

From now on let $p$ be a prime and $q\coloneqq p^e$ a $p$-power. We are interested in a special kind of affine unitals, namely affine $\SL(2,q)$-unitals. The construction of those affine unitals is due to Grundhöfer, Stroppel and Van Maldeghem \cite{slu}. They consider translations of unitals, i.\,e.\ automorphisms fixing each block through a given point (the so-called center). Of special interest are unitals of order $q$ where two points are centers of translation groups of order $q$. In the classical (Hermitian) unital of order $q$, any two such translation groups generate a group isomorphic to $\SL(2,q)$; see \cite[Main Theorem]{moufang} for further possibilities. The construction of (affine) $\SL(2,q)$-unitals is motivated by this action of $\SL(2,q)$ on the classical unital.\bigskip

Let $S\leq \SL(2,q)$ be a subgroup of order $q+1$ and let $T\leq \SL(2,q)$ be a Sylow $p$-subgroup. Recall that $T$ has order $q$ (and thus trivial intersection with $S$), that any two conjugates $T^h\coloneqq\inv hTh$, $h\in \SL(2,q)$, have trivial intersection unless they coincide and that there are $q+1$ conjugates of $T$.

Consider a collection $\mathcal D$ of subsets of $\SL(2,q)$ such that each set $D\in \mathcal D$ contains $\mathds{1}\coloneqq\mat 1001$, that $\# D=q+1$ for each $D\in \mathcal D$, and the following properties hold:

\begin{itemize}
	\item[(Q)] For each $D\in \mathcal D$, the map
	\[(D\times D)\smallsetminus \{(x,x) \mid x\in D\} \to \SL(2,q)\text{,}\quad (x,y)\mapsto x\inv y\text{,}\]
	is injective, i.\,e.\ the set $D^*\coloneqq\{x\inv y \mid x,y\in D,\ x\neq y\}$
	contains $q(q+1)$ elements.
	
	\item[(P)] The system consisting of $S\smallsetminus\{\mathds{1}\}$, all conjugates of $T\smallsetminus\{\mathds{1}\}$ and all sets $D^*$ with $D\in \mathcal D$ forms a partition of $\SL(2,q)\smallsetminus\{\mathds{1}\}$.
\end{itemize}

Set
\begin{align*}
\mathcal P \coloneqq &\ \SL(2,q)\text{,}\\
\mathcal B \coloneqq &\ \{Sg \mid g\in \SL(2,q)\} \cup \{T^h g \mid h,g\in \SL(2,q)\} \cup \{Dg \mid D\in\mathcal D, g\in\SL(2,q)\}
\end{align*}
and let the incidence relation $I\subseteq \mathcal P\times\mathcal B$ be containment.\bigskip

Then we call the incidence structure $\mathbb U_{S,\mathcal D}\coloneqq (\mathcal P,\mathcal B,I)$ an \textbf{\boldmath affine $\SL(2,q)$-unital}. Each affine $\SL(2,q)$-unital is indeed an affine unital of order $q$, see \cite[Prop. 3.15]{diss}. We call the sets $\hat D\coloneqq\{D\inv d \mid d\in D\}$, $D\in\mathcal D$, the \textbf{hats} of $\mathbb U_{S,\mathcal D}$ and the blocks $Dg$, $D\in\mathcal D$ and $g\in\SL(2,q)$, the \textbf{arcuate blocks} of $\mathbb U_{S,\mathcal D}$.\bigskip

For the construction of an affine $\SL(2,q)$-unital, we have to choose a subgroup $S\leq\SL(2,q)$ of order $q+1$ and find a set $\mathcal D$ of arcuate blocks through $\mathds{1}$ such that (Q) and (P) hold.

\goodbreak

\begin{ex}\leavevmode
	\begin{enumerate}[(a)]
		\item For each prime power $q$ we may choose $S=C$ to be cyclic and $\mathcal H$ a set of arcuate blocks through $\mathds{1}$ such that $\mathbb U_{C,\mathcal H}$ is isomorphic to the affine part of the classical unital. We call $\mathbb U_{C,\mathcal H}$ the \textbf{classical affine $\SL(2,q)$-unital}. See \cite[Example 3.1]{slu} or \cite[Section 3.2.2]{diss} for details.
		\item In \cite{slu}, Grundhöfer, Stroppel and Van Maldeghem introduce a non-classical affine $\SL(2,4)$-unital.
	\end{enumerate}	
\end{ex} 

\begin{prop}\label{gruppes}
	Let $p=2$ and let $S\leq\SL(2,q)$ be a subgroup of order $q+1$. Then $S$ is cyclic and unique up to conjugation.
\end{prop}

\begin{myproof} For $p=2$, we have $\SL(2,q)\cong \PSL(2,q)$. Using Dickson's list of subgroups of $\PSL(2,q)$ (see e.\,g.\ \cite[Hauptsatz II.8.27]{huppert}), we see that each subgroup of order $q+1$ is cyclic. From \cite[Satz II.8.5]{huppert}, we get that there is exactly one conjugacy class of cyclic subgroups of $\PSL(2,q)$ of order $q+1$.
\end{myproof}

\begin{rem}
	In \emph{\cite[Proposition 2.5]{diss}}, we give a complete list of possible subgroups $S\leq \SL(2,q)$ of order $q+1$. For $p\not\equiv 3\mod 4$, the group $S$ is cyclic. For $p\equiv 3\mod 4$, $S$ is cyclic or generalized quaternion and there is one exceptional case for $q=23$ and one for $q=47$.
\end{rem}

For each prime power $q$, we may choose a cyclic subgroup $C\leq \SL(2,q)$ of order $q+1$ as given in the following

\begin{rem}\label{sinsl}
	Let $d\in\mathbb F_q^\times$ such that $X^2-tX+d$ has no root in $\mathbb F_q$, where $t=1$ if $q$ is even and $t=0$ if $q$ is odd. Then
	\[ C\coloneqq \set{\mat ab{-db}{a+tb}}{a^2+tab+db^2=1} \]
	is a cyclic subgroup of $\SL(2,q)$ of order $q+1$. Note that $C$ is the norm $1$ group of the quadratic extension field 
	\[ \mathbb F_{q^2}\coloneqq \medset{\mat ab{-db}{a+tb}}{a,b\in\mathbb F_q} \text{.} \]
\end{rem}

We take a brief look on automorphisms of affine $\SL(2,q)$-unitals, i.\,e.\ bijections of the point set such that the block set is invariant. On any affine $\SL(2,q)$-unital $\mathbb U_{S,\mathcal D}$, right multiplication with elements of $\SL(2,q)$ obviously induces automorphisms. Let $R \coloneqq \{\rho_h \mid h\in\SL(2,q)\}\leq\Aut(\mathbb U_{S,\mathcal D})$, where $\rho_h \in R$ acts on $\mathbb U_{S,\mathcal D}$ by right multiplication with $h\in \SL(2,q)$. Every automorphism of $\SL(2,q)$ obviously induces a bijection of the point set of $\mathbb U_{S,\mathcal D}$, but it need not leave the block set invariant. Let $\mathfrak A$ denote the permutation group given by all automorphisms of $\SL(2,q)$. \bigskip

We import a useful statement from \cite{autos}:

\begin{thm}[\cite{autos}, Theorem 3.3]\label{isom}  Let $q\geq 3$ and let $\mathbb U_{S,\mathcal D}$ and $\mathbb U_{S',\mathcal D'}$ be affine $\SL(2,q)$-unitals.
	\begin{enumerate}[\em (a)]
	\item \label{isoma} Let $\psi\colon \mathbb U_{S,\mathcal D} \to \mathbb U_{S',\mathcal D'}$ be an isomorphism. Then $\psi = \alpha\rho_h$ with $\rho_h\in R$ and $\alpha\in \mathfrak A$ such that $S\cdot \alpha=S'$.
	\item \label{isomb} $\Aut(\mathbb U_{S,\mathcal D})\leq \mathfrak A_S\ltimes R$. \qed
	\end{enumerate}
\end{thm}

\begin{rem}
	The classical affine $\SL(2,q)$-unital $\mathbb U_{C,\mathcal H}$ admits the whole group $\mathfrak A_C\ltimes R$ as automorphism group \emph{(}see \emph{\cite[Proposition 4.6]{diss})}. Hence, $\mathfrak A_S\ltimes R$ is a sharp upper bound for the automorphism group of any affine $\SL(2,q)$-unital $\mathbb U_{S,\mathcal D}$ of order $q\geq 3$.
\end{rem}

\section{Three Affine \texorpdfstring{$\SL(2,8)$-Unitals}{SL(2,8)-Unitals}}

Let $q=8$ and $\mathbb F_8^\times = \langle z \rangle$, with $z^3=z+1$. The polynomial $X^2+X+1$ has no root in $\mathbb F_8$ and the Frobenius automorphism
\[\varphi\colon \mathbb F_8 \to \mathbb F_8\text{,}\quad x\mapsto x^2\text{,}\] 
has order $3$.  
Since $q=8$ is even, any subgroup $S\leq\SL(2,8)$ of order $9$ is cyclic and we may hence choose 
\[ S\coloneqq C = \set{\mat abb{a+b}}{a,b\in\mathbb F_8,\ a^2+ab+b^2=1}\text{.} \]
A generator of $C$ is given by $g\coloneqq\mat{z^2}{z^4}{z^4}{z}$. Let $f\coloneqq \mat 0110$. Then
\[ \mathfrak A_C = \Aut(\SL(2,8))_C = \langle \gamma_g\rangle \rtimes \langle \gamma_f\cdot  \varphi\rangle \cong C_9\rtimes C_6\text{,} \]
where $\varphi$ acts entrywise on a matrix and $\gamma_x$ describes conjugation with $x$. Representatives of the conjugacy classes of minimal subgroups of $\mathfrak A_C$ are
\begin{itemize}
	\item[] $F\coloneqq \langle \gamma_f \rangle \cong \langle \mat0110 \rangle \cong C_2$,
	\item[] $U\coloneqq \langle \gamma_{g^3} \rangle \cong \langle \mat1110 \rangle\cong C_3$ and
	\item[] $L\coloneqq \langle \varphi \rangle \cong C_3$.
\end{itemize}

\begin{ex}[The classical affine unital of order \unboldmath $8$]
	Let
	\begin{align*}
	H_1&\coloneqq\{ \mathds{1}, \mat {z^5}1{z^5}{z^6}, \mat {z^4}{z^2}1{z^2}, \mat 0z{z^6}{z^5}, \mat {z^3}{z^6}{z^4}{z^5}, \mat {z^3}z{z^6}0, \mat 1{z^2}1{z^6}, \mat {z^4}1{z^5}1, \mat {z^5}00{z^2}  \}\text{,}\\
	H_2&\coloneqq H_1\cdot\varphi\text{, } H_3\coloneqq H_1\cdot\varphi^2\text{,}\\
	H_4&\coloneqq\{ \mathds{1}, \mat {z^5}0{z^6}{z^2}, \mat {z}{z^6}{z^4}1, \mat 1{z^5}{z^6}{z^5}\text{,} \mat {z}{z^4}0{z^6}, \mat {z^5}{z^2}{z^4}{z^4}, \mat 0{z^2}{z^5}{z^5}, \mat 0{z^4}{z^3}{z^4}, \mat {z^2}{z^5}{z^5}{z^6}  \}\text{,}\\
	H_5&\coloneqq H_4\cdot\varphi\text{, } H_6\coloneqq H_4\cdot\varphi^2
	\end{align*}
	and $\mathcal H\coloneqq\{H_1,\ldots, H_6\}$. Then $\mathbb U_{C,\mathcal H}$ is the classical affine unital of order 8. Recall that for $H\in\mathcal H$, we denote by $\hat H$ the set of arcuate blocks $\{H\inv h\mid h\in H\}$. As indicated, $\varphi$ acts on the set of hats $\{\hat H \mid H\in\mathcal H\}$ in two orbits of length $3$. Conjugation by $g$ stabilizes each $\hat H$ and acts transitively on the blocks of each $\hat H$. Conjugation by $f$ also stabilizes each $\hat H$ but fixes exactly one block per $\hat H$. 
\end{ex}

\begin{thm}[Weihnachtsunital] Let $C\coloneqq \langle g \rangle = \langle  \left(\begin{smallmatrix} z^2&z^4\\z^4&z\end{smallmatrix}\right) \rangle$ as above and let
	\begin{align*}
	D_1&\coloneqq\{ \mathds{1}, \mat {z^5}1{z^5}{z^6}, \mat {z^4}{z^2}1{z^2}, \mat 0z{z^6}{z^2}, \mat {1}{z^4}{z^2}{z^2}, \mat {1}z{z^6}0, \mat 1{z^2}1{z^6}, \mat {z^4}1{z^5}1, \mat {z^5}00{z^2}  \}\text{,}\\
	D_2&\coloneqq D_1\cdot\varphi\text{, } D_3\coloneqq D_1\cdot\varphi^2\text{,}\\
	D_4&\coloneqq\{ \mathds{1}, \mat {z^5}0{z^6}{z^2}, \mat {z}{z^6}{z^4}1, \mat 0{z}{z^6}{z^5}, \mat {z^4}0{z^2}{z^3}, \mat {z^5}{z^2}{z^3}{z^6}, \mat 0{z^2}{z^5}{z^5}, \mat 0{z^4}{z^3}{z^4}, \mat {z^2}{z^5}{z^5}{z^6}  \}\text{,}\\
	D_5&\coloneqq D_4\cdot\varphi\text{, } D_6\coloneqq D_4\cdot\varphi^2
	\end{align*}
	and $\mathcal D\coloneqq\{D_1,\ldots, D_6\}$. Then $\mathbb{WU}\coloneqq\mathbb U_{C,\mathcal D}$ is an affine $\SL(2,8)$-unital and we call it \emph{\textbf{Weihnachtsunital}%
		\footnote{%
			The Weihnachtsunital was discovered around Christmas 2017, whence the name.}}. The stabilizer of $\mathds{1}$ in $\Aut(\mathbb{WU})$ is
	\[ \Aut(\mathbb{WU})_\mathds{1} = U \rtimes (F\times L) = \langle \gamma_{g^3} \rangle \rtimes \langle \gamma_f\cdot  \varphi \rangle \cong C_3\rtimes C_6\]
	and the full automorphism group
	\[ \Aut(\mathbb{WU}) = \Aut(\mathbb{WU})_\mathds{1} \ltimes R \]
	has index $3$ in $\Aut(\mathbb U_{C,\mathcal H})=\mathfrak A_C\ltimes R$.
\end{thm}

\begin{myproof}
	The proof is basically computation (recall Theorem \ref{isom}). Note that the given description already uses the automorphism $\varphi\in\Aut(\mathbb{WU})_\mathds{1}$. Conjugation by $f$ stabilizes each hat with exactly one fixed block per hat. Conjugation by the generator $g$ of $C$ does not induce an automorphism of $\mathbb{WU}$, but conjugation by $g^3$ yields an automorphism of~$\mathbb{WU}$ such that each hat is fixed.
\end{myproof}

Having computed the full automorphism group of $\Aut(\mathbb{WU})$, we know in particular that the Weihnachtsunital is not isomorphic to the classical affine $\SL(2,8)$-unital $\mathbb U_{C,\mathcal H}$.
Another way to see that $\mathbb{WU}$ is not isomorphic to $\mathbb U_{C,\mathcal H}$ is via O'Nan configurations. An O'Nan configuration consists of four distinct blocks meeting in six distinct points:

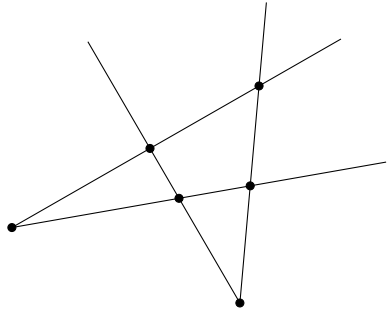
\begin{center}
	\begin{tikzpicture}
	\draw (0,0)--(10:5) (0,0)--(30:5);
	\draw (3,-1)--++(85:4) (3,-1)--++(120:4);
	\draw (0,0) coordinate (a) (10:5) coordinate (b) (30:5) coordinate (c) (3,-1) coordinate (d) ++(85:4) coordinate (e) (3,-1)++(120:4) coordinate (f);
	\draw[fill=black] (a) circle (1.5pt) (intersection of a--b and d--e) circle (1.5pt) (intersection of a--b and d--f) circle (1.5pt) (intersection of a--c and d--e) circle (1.5pt) (intersection of a--c and d--f) circle (1.5pt) (d) circle (1.5pt);
	\end{tikzpicture}
\end{center}

O'Nan observed that classical unitals do not contain such configurations (see \cite[507]{onan}).

\begin{rem}
	In $\mathbb{WU}$, there are lots of O'Nan configurations, e.\,g.\ 
	\begin{align*}
	C=\ &\{ \mathds{1}, g,g^2,g^3,g^4,g^5,g^6,g^7,g^8 \}\text{,}\\
	T\coloneqq\ &\{\mathds{1}, \mat 1101, \mat 1z01,\mat 1{z^2}01, \mat 1{z^3}01, \mat 1{z^4}01, \mat 1{z^5}01, \mat 1{z^6}01 \}\text{,}\\
	D_2\cdot\mat{z^4}110 =\ &\{ \mat {z^4}110, \mat 0{z^3}{z^4}{z^3}, \mat 1{z}01, \mat {z^2}0{z}{z^5}, \mat {z^2}1{z^2}{z^4}, \mat {z}1{z^2}{z^5}, \mat 0111, \mat {z^4}z0{z^3}, \mat 1{z^3}{z^4}0  \}\text{,}\\
	D_3\cdot\mat z{z^3}0{z^6} =\ &\{ \mat z{z^3}0{z^6}, \mat 1110, \mat {z^3}{z^4}zz, \mat 0{z^3}{z^4}{z^2}, \mat {z}1{z^2}{z^5}, \mat {z}0{z^4}{z^6}, \mat zzz{z^5}, \mat {z^3}z11, \mat 1{z^2}01 \}\text{.}
	\end{align*}
	
	\begin{center}
		\begin{tikzpicture}
		\draw[font=\small] (0,0)--(10:5) node[right,xshift=1ex]{$C$} (0,0)--(30:5) node[right,xshift=1ex,yshift=1ex]{$T$};
		\draw[font=\small] (3,-1)--++(85:4) node[above]{$D_2\cdot\mat{z^4}110$} (3,-1)--++(120:4) node[above]{$D_3\cdot\mat z{z^3}0{z^6}$};
		\draw (0,0) coordinate (a) (10:5) coordinate (b) (30:5) coordinate (c) (3,-1) coordinate (d) ++(85:4) coordinate (e) (3,-1)++(120:4) coordinate (f);
		\draw[fill=black, font=\small] (a) circle (1.5pt) node[left]{$\mathds{1}$} (intersection of a--b and d--e) circle (1.5pt) node[below right]{$g^6$} (intersection of a--b and d--f) circle (1.5pt) node[below left]{$g^3$} (intersection of a--c and d--e) circle (1.5pt) node[xshift=3ex, yshift=-1ex]{$\mat 1z01$} (intersection of a--c and d--f) circle (1.5pt) node[xshift=-4ex,yshift=1ex]{$\mat 1{z^2}01$} (d) circle (1.5pt) node[below]{$\mat {z}1{z^2}{z^5}$};
		\end{tikzpicture}
	\end{center}
\end{rem}

\goodbreak
\begin{thm}[Osterunital and Pfingstunital%
	\footnote{%
		The Osterunital and the Pfingstunital were discovered in 2018, you might guess the dates.}] Let $C\coloneqq \langle g \rangle = \langle  \left(\begin{smallmatrix} z^2&z^4\\z^4&z\end{smallmatrix}\right) \rangle$ as above. 
	\begin{enumerate}[\em (a)]
		\item Let
		\begin{align*}
		D_1&\coloneqq\{ \mathds{1}, \mat {z^5}1{z^5}{z^6}, \mat {z^4}{z^2}1{z^2}, \mat {1}z{z^6}0, \mat 0z{z^6}{z^2}, \mat {1}{z^4}{z^2}{z^2}, \mat {z^3}{z^5}{z^3}1, \mat {z^5}{z^4}{z^2}{z^4}, \mat {z^2}00{z^5}  \}\text{,}\\
		D_2&\coloneqq D_1^g\text{, } D_3\coloneqq D_1^{g^2}\text{,}\\
		D_4&\coloneqq\{ \mathds{1}, \mat {z^5}0{z^6}{z^2}, \mat {z}{z^6}{z^4}1, \mat {z^5}{z^2}{z^5}0, \mat {z^3}{z^4}{z^6}{z^5}, \mat 11{z^3}z, \mat 1z1{z^3}, \mat z{z^2}1{z^5}, \mat z0z{z^6}  \}\text{,}\\
		D_5&\coloneqq D_4^g\text{, } D_6\coloneqq D_4^{g^2}
		\end{align*}
		and $\mathcal D\coloneqq\{D_1,\ldots, D_6\}$. Then $\mathbb{OU}\coloneqq\mathbb U_{C,\mathcal D}$ is an affine $\SL(2,8)$-unital and we call it \emph{\textbf{Osterunital}}.
		
		\item Let $f=\mat0110$ as above and let
		\begin{align*}
		D_1'\coloneqq D_1\text{, } D_2' \coloneqq D_2\text{, } D_3' \coloneqq D_3\text{,}\\
		D_4'\coloneqq D_4^f\text{, } D_5'\coloneqq (D_4')^g\text{, } D_6'\coloneqq (D_4')^{g^2}
		\end{align*}
		and $\mathcal D'\coloneqq\{D_1',\ldots, D_6'\}$. Then $\mathbb{PU}\coloneqq\mathbb U_{C,\mathcal D'}$ is an affine $\SL(2,8)$-unital and we call it \emph{\textbf{Pfingstunital}}.
	\end{enumerate}
	We denote by $C$ also the automorphism group $C\coloneqq \langle \gamma_g \rangle \leq \mathfrak A_C$. The full stabilizers of~$\mathds{1}$ in $\Aut(\mathbb{OU})$ and $\Aut(\mathbb{PU})$, respectively, are
	\[ \Aut(\mathbb{OU})_\mathds{1} = \Aut(\mathbb{PU})_\mathds{1} = C \rtimes L = \langle \gamma_g \rangle \rtimes \langle \varphi \rangle \cong C_9\rtimes C_3\]
	and the full automorphism groups
	\[ \Aut(\mathbb{OU}) = \Aut(\mathbb{PU}) = (C \rtimes L) \ltimes R \]
	have index $2$ in $\Aut(\mathbb U_{C,\mathcal H})$.
\end{thm}

\begin{myproof}
	Again this is basically computation. The given description already uses the automorphism $\gamma_g$ in both $\Aut(\mathbb{OU})_\mathds{1}$ and $\Aut(\mathbb{PU})_\mathds{1}$. The Frobenius automorphism $\varphi$ acts as automorphism on $\mathbb{OU}$ as well as on $\mathbb{PU}$ in the same way as it does on $\mathbb U_{C,\mathcal H}$ and on $\mathbb{WU}$. The orbits of $\varphi$ in $\mathcal D$ are $\{D_1,D_2,D_3\}$ and $\{D_4,D_5,D_6\}$  and its orbits in $\mathcal D'$ are $\{D_1',D_2',D_3'\}$ and $\{D_4',D_5',D_6'\}$. Conjugation by $f$ induces no automorphism on neither $\mathbb{OU}$ nor $\mathbb{PU}$.
\end{myproof}

\begin{rem}
	Other than in the Weihnachtsunital, there is a difference between the action of $\Aut(\mathbb{OU})_\mathds{1}=\Aut(\mathbb{PU})_\mathds{1}\leq \mathfrak A_C$ on the set of hats of the Oster- and Pfingstunital, respectively, and its action on the set of hats of the classical affine $\SL(2,8)$-unital $\mathbb U_{C,\mathcal H}$. In $\mathbb U_{C,\mathcal H}$, conjugation by $g$ fixes every hat, while on $\mathbb{OU}$ and $\mathbb{PU}$ it acts on the set of hats in two orbits of length $3$.
\end{rem}

\begin{rem}
	As in the Weihnachtsunital, there are also many O'Nan configurations in $\mathbb{OU}$ and $\mathbb{PU}$, e.\,g.\
	\begin{align*}
	C&=\{ \mathds{1}, g,g^2,g^3,g^4,g^5,g^6,g^7,g^8 \}\text{,}\\
	D_1&=\{ \mathds{1}, \mat {z^5}1{z^5}{z^6}, \mat {z^4}{z^2}1{z^2}, \mat {1}z{z^6}0, \mat 0z{z^6}{z^2}, \mat {1}{z^4}{z^2}{z^2}, \mat {z^3}{z^5}{z^3}1, \mat {z^5}{z^4}{z^2}{z^4}, \mat {z^2}00{z^5}  \}\text{,}\\
	D_2\cdot g &=\{ g, \mat {z^4}{z^2}1{z^2}, \mat {z^3}0{z^3}{z^4}, \mat {z^5}{z^5}{z^3}{z^5}, \mat {z^6}01{z}, \mat {z^6}z{z^3}{z^6}, \mat {z^5}{z^2}{z^6}{z^5}, \mat 1{z^3}{z^4}0, \mat {z^2}{z^3}0{z^5}  \}\text{,}\\
	D_3\cdot\mat {z^5}{z^2}{z^6}{z^5} &=\{ \mat {z^5}{z^2}{z^6}{z^5}, \mat {z^2}00{z^5}, \mat z{z^2}{z^3}{z^3}, \mat {z^2}1{z^6}1, \mat {z^6}z{z^5}{z^3}, \mat {z^4}{z^3}{z^4}0, \mat z{z^4}{z^4}{z^2}, \mat {z^3}{z^2}{z^5}0, \mat {z^3}z{z^4}z \}\text{.}
	\end{align*}
	
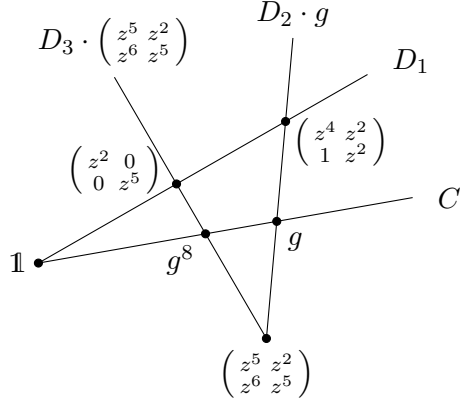
	\begin{center}
		\begin{tikzpicture}
		\draw[font=\small] (0,0)--(10:5) node[right,xshift=1ex]{$C$} (0,0)--(30:5) node[right,xshift=1ex,yshift=1ex]{$D_1$};
		\draw[font=\small] (3,-1)--++(85:4) node[above]{$D_2\cdot g$} (3,-1)--++(120:4) node[above]{$D_3\cdot\mat {z^5}{z^2}{z^6}{z^5}$};
		\draw (0,0) coordinate (a) (10:5) coordinate (b) (30:5) coordinate (c) (3,-1) coordinate (d) ++(85:4) coordinate (e) (3,-1)++(120:4) coordinate (f);
		\draw[fill=black, font=\small] (a) circle (1.5pt) node[left]{$\mathds{1}$} (intersection of a--b and d--e) circle (1.5pt) node[below right]{$g$} (intersection of a--b and d--f) circle (1.5pt) node[below left]{$g^8$} (intersection of a--c and d--e) circle (1.5pt) node[xshift=4ex, yshift=-1.5ex]{$\mat{z^4}{z^2}1{z^2}$} (intersection of a--c and d--f) circle (1.5pt) node[xshift=-4.5ex,yshift=1ex]{$\mat{z^2}00{z^5}$} (d) circle (1.5pt) node[below]{$\mat {z^5}{z^2}{z^6}{z^5}$};
		\end{tikzpicture}
	\end{center}
\end{rem}

\bigskip

Although they look quite similar, the Osterunital and the Pfingstunital are not isomorphic, as is shown in the following

\begin{prop}
	There is no isomorphism between $\mathbb{OU}$ and $\mathbb{PU}$.
\end{prop}

\begin{myproof}
	According to Theorem \ref{isom}, any isomorphism between $\mathbb{OU}$ and $\mathbb{PU}$ must be contained in $\mathfrak A_C\ltimes R$. But since the index of $\Aut(\mathbb{OU})$ in $\mathfrak A_C\ltimes R$ equals $2$ and computation shows that $D_1^f$ is no block of $\mathbb{PU}$, the statement follows.
\end{myproof}

In particular, the Oster- and Pfingstunital are two non-isomorphic affine $\SL(2,q)$-unitals with the same full automorphism group.\bigskip

\begin{rem}
	The Weihnachts-, Oster- and Pfingstunital were found by a computer search. In fact, we did an exhaustive search for affine $\SL(2,8)$-unitals, where the groups $F$, $U$ and $L$ act in the same way as on the classical affine $\SL(2,8)$-unital. Those three affine unitals were the only ones appearing through the search. See \emph{\cite[Chapter 6]{diss}} for details about the search.
\end{rem}

\section{Completion to Unitals}
Any affine unital can be completed to a unital by each of its parallelisms. In any affine $\SL(2,q)$-unital, the set of short blocks is the set of all right cosets of the $q+1$ Sylow $p$-subgroups of $\SL(2,q)$. Note that each right coset $Tg$ is a left coset $gT^g$ of a conjugate of $T$. A parallelism as in (AU5) means a partition of the set of short blocks into $q+1$ sets of $q^2-1$ pairwise non-intersecting cosets. For each prime power $q$, there are hence two obvious parallelisms, namely partitioning the set of short blocks into the sets of \emph{right} cosets or into the sets of \emph{left} cosets of the Sylow $p$-subgroups. We name those two parallelisms ``flat'' and ``natural'', respectively, and denote them by the corresponding musical signs
\[ \flat \coloneqq \{\{Tg\mid g\in\SL(2,q)\}\mid T\in\mathfrak P\}\quad \text{and}\quad \natural \coloneqq \{\{gT\mid g\in\SL(2,q)\}\mid T\in\mathfrak P\}\text{,} \]
where $\mathfrak P$ denotes the set of Sylow $p$-subgroups of $\SL(2,q)$.\bigskip

Given an affine $\SL(2,q)$-unital $\mathbb U_{S,\mathcal D}$ with parallelism $\pi$, we call the $\pi$-closure an \textbf{$\SL(2,q)$-($\pi$-)unital}. Completing $\mathbb{WU}$, $\mathbb{OU}$ and $\mathbb{PU}$ with $\flat$ and $\natural$ each, we obtain six pairwise non-isomorphic $\SL(2,q)$-unitals of order $8$. Since they are all $\flat$- or $\natural$-closures of non-classical affine $\SL(2,q)$-unitals of order $8\geq 3$, we know from \cite[Proposition 3.11 and Theorem~3.16]{autos} that their full automorphism groups fix the block $[\infty]$. Since the parallelisms $\flat$~and~$\natural$, respectively, are preserved under the action of $\mathfrak A\ltimes R$, we get
\[ \Aut(\mathbb U^\pi)= \Aut(\mathbb U^\pi)_{[\infty]} = \Aut(\mathbb U) \]
for any $\mathbb U\in \{\mathbb{WU},\mathbb{OU},\mathbb{PU}\}$ and $\pi\in\{\flat,\natural\}$.\bigskip

\begin{rem}
	In any $\SL(2,q)$-$\natural$-unital, the Sylow $p$-subgroups act (via right multiplication) as translation groups of order $q$ with centers on the block $[\infty]$. Hence, $\mathbb{WU}^\natural$, $\mathbb{OU}^\natural$ and $\mathbb{PU}^\natural$ are examples of non-classical unitals of order $q$ where the translations generate $\SL(2,q)$.
\end{rem}

\begin{rem}
	There might be more parallelisms on the short blocks of $\SL(2,8)$-unitals, leading to further closures. We already know a class of parallelisms for each odd order and one for square order \emph{(}described in \emph{\cite[Sections 2.1 and 2.2]{paratrans})} and some parallelisms for order~$4$, leading to $12$ new $\SL(2,4)$-unitals, the so-called \emph{Leonids unitals} \emph{(}see \emph{\cite[Section~2.3]{paratrans}} and \emph{\cite[Section 6.2.2]{diss})}.
\end{rem}\bigskip\bigskip

\textbf{Acknowledgment.} The author wishes to warmly thank her thesis advisor Markus J. Stroppel for his highly valuable support in each phase of this research.
	
\goodbreak	
	\printbibliography
\end{document}